\documentclass{article}

\usepackage{amsmath,amssymb,amsfonts,amsthm,latexsym,url,bm}

\newtheorem{thm}{Theorem}[section]
\newtheorem{theorem}[thm]{Theorem}

\newtheorem{cor}[thm]{Corollary}
\newtheorem{prop}[thm]{Proposition}

\newtheorem{lemma}[thm]{Lemma}

\usepackage{varioref}
\labelformat{section}{Section~#1}
\labelformat{figure}{Figure~#1}
\labelformat{table}{Table~#1}

\usepackage[numbers]{natbib}
\date{}

\begin{document}

\title{Carries, Shuffling and An Amazing Matrix}

\author{Persi Diaconis\footnote{\textit{Departments of Mathematics and Statistics,
 Stanford University}}\and%
        Jason Fulman\footnote{\textit{Department of Mathematics,
         University of Southern California; fulman@usc.edu}}}

\date{Version of June 17, 2008}

\maketitle

\begin{abstract}
The number of ``carries'' when $n$ random integers are added forms a
Markov chain \citep{Holte}. We show that this Markov chain has the
same transition matrix as the descent process when a deck of $n$ cards
is repeatedly riffle shuffled. This gives new results for the
statistics of carries and shuffling.
\end{abstract}

\section{Introduction}\label{sec1}

In a wonderful article in this monthly, John Holte \cite{Holte} found fascinating
mathematics in the usual process of ``carries'' when adding
integers. His article reminded us of the mathematics of shuffling
cards. This connection is developed below.

Consider adding two 50-digit binary numbers:
\begin{equation*}
\begin{array}{lllllllllll}
1&11111&11100&01110&01000&00001&00111&10111&00000&01111&1110\\
 &01101&11110&10111&00110&00000&10011&11011&10001&00011&11010\\
 &10111&01011&00011&10101&11110&10001&01000&11010&10101&01111\\ \hline
1&00101&01001&11010&11011&11111&00101&00100&01011&11001&01001
\end{array}
\end{equation*}
For this example, 28/50=56\% of the columns have a carry of 1. Holte shows
that if the binary digits are chosen at random, uniformly, in the limit
50\% of all the carries are zero. This holds no matter what the base. More
generally, if $n$ integers (base $b$) are produced by choosing their
digits uniformly at random in $\{0,1,\cdots,b-1\}$, the sequence of
carries $\kappa_0=0,\kappa_1,\kappa_2,\cdots$ is a Markov chain taking
values in $\{0,1,2,\cdots,n-1\}$. Holte begins by deriving the transition
matrix between successive carries $\kappa,\kappa'$.
\begin{description}
\item[(H1)] $P(i,j)=$
\begin{align*}
P(\kappa'=j|\kappa=i)&=P\left\{jb\leq i+X_1+\cdots+X_n\leq(j+1)b-1\right\}\\
&=\frac{1}{b^n}\sum_{r=0}^{j-\lfloor i/b\rfloor}(-1)^r\binom{n+1}{r}\binom{n-1-i+(j+1-r)b}{n}\notag
\end{align*}
Here, $0\leq i,j \leq n-1$ and $X_1,X_2,\cdots,X_n$ are independent and
uniformly distributed on $\{0,1,\cdots,b-1\}$.

\item[(H2)] When $b=2$, for any $n$, the transition matrix is
\begin{equation*}
P(i,j)=\frac{1}{2^n}\cdot\binom{n+1}{2j-i+1}\qquad 0\leq i,j\leq n-1.
\end{equation*}

\item[(H3)] For $n=3$, for all $b$
\begin{equation*}
P(i,j)=\frac{1}{6b^2}\begin{pmatrix}
b^2+3b+2&4b^2-4&b^2-3b+2\\
b^2-1&4b^2+2&b^2-1\\
b^2-3b+2&4b^2-4&b^2+3b+2
\end{pmatrix}.
\end{equation*}
\end{description}
These are the ``amazing matrices'' of Holte's title. Among many
things, Holte shows
\begin{description}
\item[(H4)] The matrix $P(i,j)$ of \textbf{(H1)} has stationary vector
  $\pi_n(j)$ (left eigenvector with eigenvalue 1) independent of the
  base $b$:
\begin{equation*}
\pi_n(j)=\frac{A(n,j)}{n!}
\end{equation*}
with $A(n,j)$ the Eulerian number. This may be defined as

\item[(H4$'$)] $A(n,j)$ is the number of permutations in the symmetric
  group $S_n$ with $j$-descents. Recall that $\sigma\in S_n$ has a
  descent at $i$ if $\sigma(i+1)<\sigma(i)$. So $\bm{5}\,1\,\bm{3}\,2\,4$ has
  two descents.

\item[(H4$''$)] $A(n,j)$ is the coefficient of $x^{j+1}$ in the polynomial
  $p_n(x)$ where
\begin{equation*}
\sum_{i=0}^\infty i^nx^i=\frac{p_n(x)}{(1-x)^{n+1}}.
\end{equation*}

\item[(H4$'''$)] $A(n,j)=\sum_{\ell=0}^j(-1)^\ell\binom{n+1}{\ell}(j+1-\ell)^n$.
\end{description}
Definition \textbf{(H4$'$)} is most relevant to the present paper.
\textbf{(H4$''$)} is equivalent to Worpitzky's identity. It has many
proofs and appearances, e.g., to juggling sequences \cite{BEGW}. Finally,
\textbf{(H4$'''$)} goes back to Euler. An elementary development of these
ideas is in \cite{Com}.

When $n=2, A(2,0)=A(2,1)=1$, thus $\pi_2(0)=\pi_2(1)=1/2$ is the
limiting frequency of carries when two long integers are added. When
$n=3,\, A(3,0)=1,\, A(3,1)=4,\, A(3,2)=1$, giving $\pi_3(0)=1/6,\,
\pi_3(1)=2/3,\, \pi_3(2)=1/6$.

Holte further shows
\begin{description}
\item[(H5)] The matrix $P(i,j)$ of \textbf{(H1)} has eigenvalues $1,1/b,
  1/b^2,\cdots,1/b^{n-1}$ with explicitly computable eigenvectors
  independent of $b$.
\item[(H6)] Let $P_b$ denote the matrix in \textbf{(H1)}. Then for all
  real $a,b$
\begin{equation*}
P_aP_b=P_{ab}.
\end{equation*}
\end{description}

When we saw properties \textbf{(H4)}, \textbf{(H5)}, \textbf{(H6)}, we
hollered ``Wait, this is all about shuffling cards!'' Knowledgeable
readers may well think, ``For these two guys, everything is about
shuffling cards.'' While there is some truth to these thoughts, we justify
our claim in the next section. Following this we show how the connection
between carries and shuffling contributes to each subject. The rate of
convergence of the Markov chain \textbf{(H1)} to the stationary
distribution $\pi_n$ is given in \ref{sec4}: the argument shows that the
matrix $P$ is totally positive of order 2. Finally, we show how the same
matrix occurs in taking sections of generating functions \cite{BW},
discuss carries for multiplication, and describe another ``amazing
matrix''.

Our developments do not exhaust the material in Holte's article, which we
enthusiastically recommend. A ``higher math'' perspective on arithmetic
carries as cocycles \cite{Isak} suggests many further projects. We have
tried to keep the presentation elementary, and mention the (more
technical) companion paper \cite{DF} which analyzes the carries chain
using symmetric function theory and gives analogs of our main results for
other Coxeter groups.

\section{Shuffling Cards}\label{sec2}

How many times should a deck of $n$ cards be riffle shuffled to
thoroughly mix it? For an introduction to this subject, see
\cite{AD,Mann}. The main theoretical developments are in \cite{BayerD,
  DMP} with further developments in \cite{FI,F}. A survey of the many
connections and developments is in \cite{Dia}. The basic shuffling
mechanism was suggested by \cite{GSR}. It gives a realistic
mathematical model for the usual method of riffle shuffling $n$ cards:
\begin{itemize}
\item Cut off $C$ cards with probability $\binom{n}{C}/2^n,\ 0\leq
  C\leq n$.
\item Shuffle the two parts of the deck according to the following
  rule: if at some stage there are $A$ cards in one part and $B$ cards
  in the other part, drop the next card from the bottom of the first
  part with probability $A/(A+B)$ and from the bottom of the second
  part with probability $B/(A+B)$.
\item Continue until all cards are dropped.
\begin{equation*}
A\Big\{\stackrel{=}{\equiv}\ \ \equiv\Big\}B
\end{equation*}
\end{itemize}
Let $Q(\sigma)$ be the probability of generating the permutation
$\sigma$ after one shuffle, starting from the identity. Repeated shuffling
is modeled by convolution:
\begin{equation}
Q^2(\sigma)=\sum_\eta Q(\eta)Q(\sigma\eta^{-1}),\quad Q^h(\sigma)=\sum Q^{h-1}(\eta)Q(\sigma\eta^{-1}).
\label{7}
\end{equation}
Thus to be at $\sigma$ after two shuffles, the first shuffle goes to
some permutation $\eta$ and the second must be to
$\sigma\eta^{-1}$. The uniform distribution is
$U(\sigma)=1/n!$. Standard theory shows that
\begin{equation}
Q^h(\sigma)\to U(\sigma)\quad\text{as }h\to\infty.
\label{8}
\end{equation}
The references above give useful rates for the convergence in \eqref{8}
showing that it takes $h=3/2\log_2n+c$ to get $2^{-c}$ close to
random. When $n=52$, this becomes $h\doteq 7$ shuffles.

To explain the connection with carries, it is useful to have a second
description of shuffling. Consider dropping $n$ points uniformly at random
into $[0,1\}$. Label these points in order $x_{(1)}\leq x_{(2)}\cdots\leq
x_{(n)}$. The Bakers transformation $x\mapsto 2x$ (mod 1) maps $[0,1]$
into itself and permutes the points. Let $\sigma$ be the induced
permutation. As shown in \cite{BayerD}, the chance of $\sigma$ is exactly
$Q(\sigma)$. A natural generalization of this shuffling scheme to
``$b$-shuffles'' is induced from $x\mapsto bx$ (mod 1) with $b$ fixed in
$\{1,2,3,\cdots\}$. Thus ordinary riffle shuffles are 2-shuffles and a
3-shuffle results from dividing the deck into three piles and dropping
cards sequentially from the bottom of each pile with probability
proportional to packet size.

Let $Q_b(\sigma)$ be the probability of $\sigma$ after a
$b$-shuffle. From this geometric description,
\begin{equation}
Q_a*Q_b=Q_{ab}.
\label{extra8}
\end{equation}
The Gilbert--Shannon--Reeds measure is $Q_2$ in this notation and we see
that $Q_2^h=Q_{2^h}$. Thus to study repeated shuffles, we need only
understand a single $b$-shuffle. A main result of \cite{BayerD} is a
simple formula:
\begin{equation}
Q_b(\sigma)=\frac{\binom{n+b-r}{n}}{b^n}.
\label{9}
\end{equation}
Here $r=r(\sigma)=1+\#\{\text{descents in }(\sigma^{-1})\}$.

In addition to the similarities between \textbf{(H6)} and \eqref{extra8},
\cite{BayerD} and \cite{Han} proved that the eigenvalues of the Markov
chain induced by $Q_b$ are $1,1/b,1/b^2,\cdots,1/b^{n-1}$. This and the
appearance of descents convinced us that there must be an intimate
connection between carries and shuffling. The main result of this article
makes this precise.
\begin{theorem}
The number of descents in successive $b$-shuffles of $n$ cards forms a
Markov chain on $\{0,1,\cdots,n-1\}$ with transition matrix $P(i,j)$ of
\textbf{(H1)}. \label{th1}
\end{theorem}

\section{Bijective Methods} \label{sec3}

First we describe some notation to be used throughout. The number of
descents of a permutation $\tau$ is denoted by $d(\tau)$. Label the
columns of the $n$ numbers to be added mod $b$ by $C_1,C_2,C_3,\cdots$
where $C_1$ is the right-most column.

The main purpose of this section is to give a bijective proof of the
following theorem, which implies Theorem \ref{th1} from the introduction.

\begin{theorem} \label{main} Let $\kappa_j$ denote the amount carried from
column $j$ to column $j+1$ when $n$ length $m$ numbers are added mod $b$.
Let $\tau_j$ be the permutation obtained after the iteration of $j$
$b$-shuffles, started at the identity. Then
\begin{equation*}
\mathbb{P}(\kappa_1=i_1,\cdots,\kappa_m=i_m) = \mathbb{P}
(d(\tau_1)=i_1,\cdots,d(\tau_m)=i_m)
\end{equation*}
for all values of $i_1,\cdots,i_m$.
\end{theorem}

In preparation for the proof, some notation and lemmas will be needed.

\begin{lemma} \label{clear} Let $\kappa(C_j \cdots C_1)$ denote the amount
carried from column $j$ to column $j+1$ when the corresponding $j$-tuples
are added (adding consecutive $j$-tuples one at a time rather than adding
a column at a time). Then $\kappa(C_j \cdots C_1) = \kappa_j$.
\end{lemma}

\begin{proof}
This is clear since in calculating the carry to column $j+1$
it is irrelevant how one adds the numbers in the preceding columns.
\end{proof}

Given a length $n$ list of $j$-tuples of numbers mod $b$, one says that
the list has a descent at position $i$ if the $i+1$st $j$-tuple is smaller
than the $i$th $j$-tuple. For example the following
$3$-tuples of mod $3$ numbers:
\begin{equation*}
\begin{array}{c c c}
0 & 1 & 2 \\1 & 0 & 1 \\ 2 & 2 & 0 \\ 1 & 0 & 1 \\ 0 & 2 & 0 \\ 2 & 1 & 1
\end{array}
\end{equation*}
has a descent at position $3$ since $220$ is greater
than $101$, and a descent at position $4$ since $101$ is greater than
$020$.

Given a length $n$ list of $j$-tuples of numbers mod $b$, one says that
the list has a carry at position $i$ if the addition of the $i+1$st
$j$-tuple on the list to the sum of the first $i$ $j$-tuples increases the
amount that would be carried to the $j+1$st column (it might seem more
natural to say that the carry is at position $i+1$, but our convention
will be useful). For example the following $3$-tuples of mod $3$ numbers:
\begin{equation*}
\begin{array}{c c c}
0 & 1 & 2 \\ 0 & 1 & 2 \\ 1 & 1 & 2 \\ 1 & 1 & 1 \\
2 & 1 & 2 \\ 1 & 2 & 1
\end{array}
\end{equation*}
has a carry at positions 3 and 4. Indeed $(0,1,2)+(0,1,2)=(1,0,1)$ which
doesn't create a carry. Adding $(1,1,2)$ gives $(2,2,0)$ which still
doesn't create a carry. Adding $(1,1,1)$ gives $(1,0,1)$ with a carry, so
there is a carry at position 3. Adding $(2,1,2)$ gives $(0,2,0)$ with a
carry, so there is a carry at position 4. Finally adding $(1,2,1)$ gives
$(2,1,1)$, which doesn't create a carry.

For what follows we use a bijection, which we call the bar map, on sets of
$j$ column vectors having length $n$ and entries in $0,1,\cdots,b-1$.
Given $C_j \cdots C_1$, then $\overline{C_j \cdots C_1}$ is defined as
follows: the $i$th $j$-tuple of $\overline{C_j \cdots C_1}$ consists of
the right-most $j$ coordinates of the mod $b$ sum of the first $i$
$j$-tuples of $C_j \cdots C_1$. For example,
\begin{equation*} C_3 C_2 C_1 =
\begin{array}{c c c}
0 & 1 & 2 \\ 0 & 1 & 2 \\ 1 & 1 & 2 \\ 1 & 1 & 1 \\2 & 1 & 2 \\ 1 & 2& 1
\end{array} \mapsto, \overline {C_3 C_2 C_1} =\begin{array}{c c c}
0 & 1 & 2 \\ 1 & 0 & 1 \\ 2 & 2 & 0 \\
1 & 0 & 1 \\ 0 & 2 & 0 \\ 2 & 1 & 1
\end{array}.
\end{equation*}
Indeed $012+012=101$ giving the second line of $\overline{C_3 C_2
  C_1}$. Then $101+112=220$ giving the third line, and $220+111=101$
(retaining only the last 3 coordinates), giving the fourth line, etc.
One can easily invert the bar map, so it is a bijection.

The following lemma is immediate from these definitions.

\begin{lemma} \label{descar} $\overline{C_j \cdots C_1}$ has a descent at
position $i$ if and only if $C_j \cdots C_1$ has a carry at position $i$.
\end{lemma}

Given a length $n$ collection of $j$-tuples of numbers mod $b$, we define
an associated permutation $\pi$ by labeling the $j$-tuples from
lexicographically  smallest to largest (considering the higher up
$j$-tuple to be smaller in case of ties). For example with $n=6, j=2,
b=3$, one would have
\begin{equation*}
\pi \left(
\begin{array}{c c}
1 & 2 \\ 2 & 1 \\ 1 & 0 \\ 0 & 1 \\0 & 0 \\ 2 & 1
\end{array} \right) =
\begin{array}{c}
4 \\ 5 \\ 3 \\ 2 \\ 1 \\ 6
\end{array},
\end{equation*}
since
$(0,0)$ is the smallest, followed by $(0,1)$, $(1,0)$, $(1,2)$, then the
uppermost copy of $(2,1)$ and finally the lowermost copy of $(2,1)$. Note
that we use the standard convention for writing permutations, i.e. $1
\mapsto 4$, $2 \mapsto 5$, etc. We mention that this construction appears
in the theory of inverse riffle shuffling \cite{BayerD}.

\begin{lemma} \label{samedes} $\overline{C_j \cdots C_1}$ has a descent
at position $i$ if and only if the associated permutation
$\pi(\overline{C_j \cdots C_1})$ has a descent at position $i$.
\end{lemma}

\begin{proof} This is immediate from the definition of $\pi$.
\end{proof}

To proceed define a second bijection, called the star map, on sets of $j$
column vectors having length $n$ and entries in $0,1,\cdots,b-1$. This
sends column vectors $A_j \cdots A_1$ to $(A_j \cdots A_1)^*$ defined as
follows. The right-most column of $(A_j \cdots A_1)^*$ is $A_1$. The
second column in $(A_j \cdots A_1)^*$ is obtained by putting the entries
of $A_2$ in the order specified by the permutation corresponding to
right-most column of $(A_j \cdots A_1)^*$ (which is $A_1$). Then the third
column in $(A_j \cdots A_1)^*$ is obtained by putting the entries of $A_3$
in the order specified by the permutation corresponding to the two
right-most columns of $(A_j \cdots A_1)^*$, and so on.

For example,
\begin{equation*} A_3 A_2 A_1 =
\begin{array}{c c c}
1 & 2 & 2 \\ 1 & 2 & 1 \\ 2 & 0 & 0 \\
0 & 0 & 1 \\ 2 & 1 & 0 \\ 0 & 1 & 1
\end{array} \mapsto
(A_3 A_2A_1)^* =
\begin{array}{c c c}
0 & 1 & 2 \\ 1 & 0 & 1 \\ 2 & 2 & 0 \\
1 & 0 & 1 \\ 0 & 2 & 0 \\ 2 & 1 & 1
\end{array}
\end{equation*}

Indeed, the right-most column of $(A_3 A_2 A_1)^*$ is $A_1$. The second
column of $(A_3 A_2 A_1)^*$ is obtained by taking the entries of $A_2$
(namely $2,2,0,0,1,1$) and putting the $2$ next to the smallest element of
$A_1$ (so the highest 0), then the second $2$ next to the the 2nd smallest
element (so the second 0), then the $0$ next to the 3rd smallest element
(so the highest 1), then the second $0$ next to the 4th smallest element
(so the second 1), then the $1$ next to the 5th smallest element (so the
third 1), and finally the second $1$ next to the 6th smallest element (so
the only 2), giving
\begin{equation*}
\begin{array}{c c}
1 & 2 \\ 0 & 1 \\ 2 & 0 \\
0 & 1 \\ 2 & 0 \\ 1 & 1
\end{array}.
\end{equation*}
Then the third column from of $(A_3 A_2 A_1)^*$ is obtained by taking
the entries of $A_3$ (namely $1,1,2,0,2,0$) and putting the $1$ next
to the smallest pair (so the highest $(0,1)$), then putting the second
$1$ next to the 2nd smallest pair (so the second $(0,1)$), then the
$2$ next to the third smallest pair $(1,1)$, then the $0$ next to the
fourth smallest pair $(1,2)$, then the second $2$ next to the fifth
smallest pair (the highest $(2,0)$), and finally the second $0$ next
to the sixth smallest pair (the second $(2,0)$).

The star map is straightforward to invert (we leave this as an exercise to
the reader), so it is a bijection.

The crucial property of the star map is given by the following lemma, the
$j=2$ case of which is essentially equivalent to the ``$A^B \& B$''
formula in Section 9.4 of \cite{Mann}.

\begin{lemma} \label{starkey}
\begin{equation*}
\pi(A_j) \cdots \pi(A_1) = \pi[(A_j \cdots A_1)^*] ,
\end{equation*}
where the product on the left is the usual multiplication of permutations.
\end{lemma}

As an illustration,
\begin{equation*}  A_3 A_2 A_1 =
\begin{array}{c c c}
1 & 2 & 2 \\ 1 & 2 & 1 \\ 2 & 0 & 0 \\
0 & 0 & 1 \\ 2 & 1 & 0 \\ 0 & 1 & 1
\end{array}
\end{equation*}
yields the permutations
\begin{equation*}
\begin{array}{c c c}
\pi(A_3) & \pi(A_2) & \pi(A_1) \\ 3 & 5 & 6 \\ 4 & 6 & 3
\\ 5 & 1 & 1 \\ 1 & 2 & 4 \\ 6 & 3 & 2 \\ 2 & 4 & 5
\end{array}.
\end{equation*}
Also as calculated above,
\begin{equation*}
(A_3 A_2 A_1)^* =
\begin{array}{c c c}
0 & 1 & 2 \\ 1 & 0 & 1 \\ 2 & 2 & 0 \\
1 & 0 & 1 \\ 0 & 2 & 0 \\ 2 & 1 & 1
\end{array}
\end{equation*}
which yields the permutations
\begin{equation*} \begin{array}{c c c}
\pi[(A_3A_2A_1)^*] & \pi[(A_2A_1)^*] & \pi[(A_1)^*] \\ 1 & 4 & 6 \\
3 & 1 & 3 \\ 6 & 5 & 1 \\ 4 & 2 & 4 \\ 2 & 6 & 2 \\ 5 & 3 & 5
\end{array}.
\end{equation*}
$\pi(A_1^*)=\pi(A_1)$, $\pi[(A_2A_1)^*]=\pi(A_2) \pi(A_1)$, and
$\pi[(A_3A_2A_1)^*]=\pi(A_3) \pi(A_2) \pi(A_1)$, and Lemma \ref{starkey}
gives that this happens in general.

\begin{proof}[Proof of Lemma \ref{starkey}.] This is clear for $j=1$, so
  consider $j=2$. Then the claim is perhaps easiest to see using the
  theory of inverse riffle shuffles. Namely given a column of $n$
  numbers mod $b$, mark cards $1,\cdots,n$ with these numbers, then
  bring the cards labeled $0$ to the top (cards higher up remaining
  higher up), then bring the cards labeled $1$ just beneath them, and
  so on. For instance,
\begin{equation*}
\begin{array}{c}
2 \\ 1 \\ 0 \\ 1 \\ 0 \\ 1
\end{array} \mapsto
\begin{array}{c}
3 \\ 5 \\ 2 \\ 4 \\ 6 \\ 1
\end{array}.
\end{equation*}
Note that (in the notation of the example) this is $\pi(A_1)^{-1}$. Now
repeat this process, using the column
\begin{equation*} \begin{array}{c}
2 \\ 2 \\ 0 \\ 0 \\ 1 \\ 1
\end{array},
\end{equation*}
to label the cards, placing the labels just to the left of the digit
already on each card. A moment's thought shows that this is equivalent
to a single process in which one labels the cards with pairs from
$(A_2 A_1)^*$. Thus $\pi[(A_2A_1)^*]^{-1}=\pi(A_1)^{-1}
\pi(A_2)^{-1}$, so that $\pi[(A_2A_1)^*] = \pi(A_2) \pi(A_1)$. The
reader desiring further discussion for the case of two columns is
referred to Section 9.4 of the expository paper \cite{Mann}. The argument
for $j \geq 3$ is identical: just use the observation that iterating
the procedure three times is equivalent to a single process in which
one labels the cards with triples from $(A_3 A_2 A_1)^*$.
\end{proof}

With the above preparations in hand, Theorem \ref{main} can
be proved.

\begin{proof}[Proof of Theorem \ref{main}.] To begin, note that
\begin{eqnarray*}
\kappa_1=i_1, \cdots, \kappa_m=i_m & \leftrightarrow & \kappa(C_j \cdots
C_1)=i_j \ (1 \leq j \leq m) \\
& \leftrightarrow & d( \overline{C_j \cdots C_l}) = i_j \ (1 \leq j \leq
m) \\
& \leftrightarrow & d( \pi(\overline{C_j \cdots C_l})) = i_j \ (1 \leq j
\leq m).
\end{eqnarray*}
The first step used Lemma \ref{clear}, the second step used Lemma
\ref{descar} and the third step used Lemma \ref{samedes}.

Let $A_m \cdots A_1 = (\overline{C_m \cdots C_1})^{-*}$. Then $A_j \cdots
A_1 = (\overline{C_j \cdots C_1})^{-*}$ for all $1 \leq j \leq m$, and
Lemma \ref{starkey} implies that
\begin{equation*}
d [ \pi(A_j) \cdots \pi(A_1)] = d(\pi[(A_j \cdots A_1)^*]) = d
[\pi(\overline{C_j \cdots C_1})] = i_j
\end{equation*}
for all $1 \leq j \leq m$. Now note that if $C_m \cdots C_1$ are chosen
i.i.d. with entries uniform in $0,1,\cdots,b-1$, then the same is true of
$A_m \cdots A_1$ since the bar and star maps are both bijections. Note
that each $\pi(A_i)$ has the distribution of a permutation after a
b-shuffle, so one may take $\tau_j$ to be the product $\pi(A_j) \cdots
\pi(A_1)$, and the theorem is proved.
\end{proof}

{\it Remark and example:} The above construction may appear complicated,
but we mention that the star map (though useful in the proof) is not
needed in order to go from the columns of numbers being added to the
$\tau$'s. Indeed, from the proof of Theorem \ref{main} one sees that the
$\tau_j$'s can be defined by $\tau_j=\pi(\overline{C_j \cdots C_1})$. Thus
in the running example,

\begin{equation*}
C_3C_2C_1 =
\begin{array}{c c c}
& & \\ 0 & 1 & 2 \\ 0 & 1 & 2 \\ 1 & 1 & 2
\\ 1 & 1 & 1 \\ 2 & 1 & 2 \\ 1 & 2 & 1
\end{array} \mapsto \overline{C_3C_2C_1} =
\begin{array}{c c c}
& & \\
 0 & 1 & 2 \\ 1 & 0 & 1 \\ 2 & 2 & 0 \\
1 & 0 & 1 \\ 0 & 2 & 0 \\ 2 & 1 & 1
\end{array} \mapsto
\begin{array}{c c c}
\tau_3 & \tau_2 &\tau_1 \\
 1 & 4 & 6 \\ 3 & 1 & 3 \\ 6 & 5 & 1 \\
4 & 2 & 4 \\ 2 & 6 & 2 \\ 5 & 3 & 5
\end{array}.
\end{equation*}
Observe that $\kappa_1=3$, $\kappa_2=3$, $\kappa_3=2$, and that
$d(\tau_1)=3$, $d(\tau_2)=3$, $d(\tau_3)=2$ as claimed.

As a corollary of Theorem \ref{main}, we deduce that the descent process
after riffle shuffles is Markov (usually, a function of a Markov chain is
not Markov).

\begin{cor} Let a Markov chain on the symmetric group begin at the
  identity and proceed by successive independent $b$-shuffles. Then
  $d(\pi)$, the number of descents, forms a Markov chain.
\end{cor}

\begin{proof} This follows from Theorem \ref{main} and the fact that
  the carries process is Markov.
\end{proof}

\section{Applications to the Carries Process}\label{sec4}

As in previous sections, let $\kappa_j$ be the amount carried from column
$j$ to column $j+1$ when $n$ length-$m$ numbers are added mod $b$. Suppose
throughout this section that the ``digits'' of these numbers are chosen
uniformly and independently in $\{0,1,\cdots,b-1\}$.
\begin{theorem}
  For $1\leq j\leq m$, the expected value of $\kappa_j$ is
  $\mu_j= \frac{n-1}{2} \left( 1- \frac{1}{b^j} \right)$. The variance
  of $\kappa_j$ is $\sigma_j^2= \frac{n+1}{12} \left(1- \frac{1}{b^{2j}}
  \right)$. Normalized by its mean and variance, for large $n$, $\kappa_j$ has a
  limiting standard normal distribution. \label{th4.1}
\end{theorem}

\begin{proof}
From Lemma \ref{descar} of \ref{sec3}, $\kappa_j$ is distributed exactly
like the number of descents among the $n$ rows of the right-most $j$
digits of the random array. The distribution of these descents is studied
in \cite{BDFb} where they are shown to be a 2-dependent process with the
required mean and variance. The central limit theorem for 2-dependent
processes is classical \cite{And}.
\end{proof}
{\it Remarks:} \begin{enumerate}

\item Note that $\mu_j,\sigma_j^2$ are increasing to their limiting value
$\frac{n-1}{2}, \frac{n+1}{12}$ as $j$ increases.

\item Let $S_m=\kappa_1+\kappa_2+\cdots+\kappa_m$ be the total number of
carries. By linearity of expectation and Theorem \ref{th4.1}, this has
mean
\begin{equation*}
\bar\mu_m=\frac{n-1}{2}\left(m-\frac{1}{b-1}\left(1-\frac{1}{b^m}\right)\right).
\end{equation*}
When $n=2$, this was shown by Knuth \citep[p.~278]{Knuth}. He also finds
the variance of $S_m$ when $n=2$. For fixed $n$ and $b$, the central limit
theorem for finite state space Markov chains \cite{Bill} shows that $S_m$,
normalized by its mean and variance, has a standard normal limiting
distribution.

\item The fine properties of the number of carries within a column is
studied in \cite{BDFa} where it is shown to be a determinantal point
process.
\end{enumerate}

As shown above, the carries process $\kappa_j, 0\leq j\leq m$ (with
$\kappa_0=0$) is a Markov chain which has limiting stationary distribution
$\pi(j)=A(n,j)/n!$. To study the rate of convergence to the limit we first
prove a new property of the amazing matrix $P(i,j)$ of \textbf{(H1)}.
Recall that a matrix is totally positive of order two $(TP_2)$ if all the
$2\times 2$ minors are non-negative.
\begin{lemma} \label{TP2}
For every $n$ and $b$, the matrix $P(i,j)$ of \textbf{(H1)} is $TP_2$.
\end{lemma}
\begin{proof}
As noted on p. 140 of \cite{Holte},
\begin{equation*}
P(i,j)=\frac{1}{b^n} \left[x^{(j+1)b-i-1}\right]\left(\frac{1-x^b}{1-x}\right)^{n+1}
\end{equation*}
where $[x^i]f(x)$ denotes the coefficient of $x^i$ in a polynomial $f(x)$.
Thus the transpose of $P$ is a submatrix of the matrix with $(i,j)$
coordinates $[x^{i-j}][(1-x^b)/(1-x)]^{n+1}$. Since the product of $TP_2$
matrices is $TP_2$, it is enough to treat the case $n=0$. Now, the matrix
is a lower triangular, $n\times n$ matrix with ones down the diagonal,
ones on the next lowest $b-1$ diagonals and zeros elsewhere. For example,
when $n=6,b=3$ the relevant matrix is
\begin{equation*}
\begin{array}{c c c c c c}
1&0&0&0&0&0\\
1&1&0&0&0&0\\
1&1&1&0&0&0\\
0&1&1&1&0&0\\
0&0&1&1&1&0\\
0&0&0&1&1&1
\end{array}.
\end{equation*}
By inspection, 13 of the 16 possible $2\times 2$ matrices can occur as
minors. The missing ones are
\begin{equation*}
\begin{array}{c c c}
01&01&11\\
10&11&10
\end{array},
\end{equation*}
these being the only ones with negative determinants.
\end{proof}

{\it Remark:} When $b=2$, the original $P(i,j)= 2^{-n}
\binom{n+1}{2j-i+1}$ is totally positive ($TP_\infty$). Indeed,
$P(i,j)=2^{-n}[x^{2j-i+1}](1+x)^{n+1}$. Let $i'=i+1,\,j'=j+1$. This
becomes $2^{-n}[x^{2j'-i'}](1+x)^{n+1}$. Each minor of this is a subminor
of $2^{-n}[x^{j'-i'}](1+x)^{n+1}$. This is totally positive by the
classification of Polya frequency sequences due to Schoenberg and Edrei
(\cite{Karl}, Chap. 8).

\vspace{.15in}

Consider the basic transition matrix $P(i,j)$ for general $b,n$. This has
stationary distribution $\pi(j),\,0\leq j\leq n$, given in \textbf{(H4)}.
The carries Markov chain starts at $0$ and the right-most carries tend to
be smaller. This is seen in Theorem \ref{th4.1}. It is natural to ask how
far over one must go so that the carries process is stationary. If
$P^r(0,j)$ is the chance of carry $j$ after $r$ steps, we measure the
approach to stationarity by separation
\begin{equation*}
\text{sep}(r)=\max_j \left[ 1-\frac{P^r(0,j)}{\pi(j)} \right].
\end{equation*}
Thus $0\leq \text{sep}(r)\leq 1$ and sep($r$) is small provided
$P^r(0,j)$ is close to $\pi(j)$ for all $j$. See \cite{AD} or \cite{Dia} for
further properties of separation. The following theorem shows that
convergence requires $r=2\log_bn$.
\begin{theorem} \label{sharpsep}
For any $b\geq 2, n\geq 2$, the transition matrix $P(i,j)$ of
\textbf{(H1)} satisfies
\begin{enumerate}
\item For all $r\geq0$, the separation distance $\text{sep}(r)$ of the carries
chain after $r$ steps (started at $0$) is attained at the state $j=n-1$.
\item For $r=2\log_b(n)+\log_b(c)$,
\begin{equation*}
\text{sep}(r)\to 1-e^{-\frac{1}{2c}}
\end{equation*}
if $c>0$ is fixed and $n\to\infty$.
\end{enumerate}
\end{theorem}
\begin{proof} By Lemma \ref{TP2}, the matrix $P(i,j)$ is $TP_2$. Thus
the matrix $P^*(i,j):=[P(j,i)\pi(j)]/\pi(i)$ is also $TP_2$, since every
$2\times 2$ minor of $P^*$ is a positive multiple of a $2\times 2$ minor
of $P$. Now consider the function $f_r(i)=P^r(0,i)/\pi(i)$. We claim that
$P^*f_r=f_{r+1}$. Indeed,
\begin{align*}
[P^*f_r](i) &=\sum_j P^*(i,j)f_r(j)\\
&=\sum_j P^*(i,j)\frac{P^r(0,j)}{\pi(j)}\\
&=\sum_j \frac{P(j,i)\pi(j)}{\pi(i)}\ \frac{P^r(0,j)}{\pi(j)}\\
&=\frac{P^{r+1}(0,i)}{\pi(i)}.
\end{align*}

Now the ``variation-diminishing property'' (p.22 of \cite{Karl}) gives
that if $f$ is monotone and $P^*$ is $TP_2$, then $P^*f$ is monotone.
Since $f_0$ is monotone (the walk is started at $0$), it follows that
$f_r$ is monotone, i.e., that the separation distance $s(r)$ is attained
at the state $n-1$.

For the second assertion, note that by the relation between riffle
shuffling and the carries chain in Theorem \ref{main}, $P^r(0,n-1)$ is
equal to the chance of being at the unique permutation with $n-1$ descents
after $r$ iterations of a $b$-shuffle; by \cite{BayerD} this is
$b^{-rn}\binom{b^r}{n}$. Thus
\begin{align*}
\text{sep}(r) &=1-\frac{P^r(0,n-1)}{\pi(n-1)}\\
&=1-\prod_{i=1}^{n-1}\left(1-\frac{i}{b^r}\right)\\
&=1-\exp\left(\sum_{i=1}^{n-1}\log\left(1-\frac{i}{b^r}\right)\right).
\end{align*}
Letting $b^r=cn^2$ with $c>0$ fixed, this becomes
\begin{equation*}
1-\exp\left(-
\sum_{i=1}^n\frac{i}{cn^2}+O\left(\frac{i^2}{n^4}\right)\right)
\sim1-e^{-\frac{1}{2c}},
\end{equation*}
as $n\to\infty$.
\end{proof}

{\it Remark:} It is known \cite{BayerD} that it takes $r=2\log_bn$
$b$-shuffles to make separation distance small on the symmetric group. Via
Theorem \ref{main}, this shows $2\log_bn$ steps suffice for the carries
process. Of course, fewer steps might suffice but Theorem \ref{sharpsep}
shows the result is sharp for large $n$. In mild contrast, it is known
\cite{Al,BayerD} that $(3/2)\log_2n$ ``ordinary'' ($b=2$) riffle shuffles
are necessary and suffice for total variation convergence. We can show
that for $b=2$, $\log_2n$ carry steps suffice for binary addition. Our
argument uses the monotonicity proved above, the first eigenvector from
\cite{Holte}, and Proposition 2.1 of \cite{DKS}; for a second argument,
using symmetric functions, see \cite{DF}. We do not know that this upper
bound is sharp; the best total variation lower bound we have is
$(1/2)\log_2n$.

\section{Three Related Topics}\label{sec5}

The ``amazing matrix'' turns up in different contexts (sections of
generating functions) in the work of Brenti--Welker \cite{BW}. There is an
analog of carries for multiplication which has interesting structure.
Finally, there are quite different amazing matrices having many of the
same properties as Holte's. These three topics are briefly developed in
this section.

\subsection{Sections of generating functions}

Some natural sequences $a_k,\,0\leq k<\infty$ have generating functions:
\begin{equation}
\sum_{k=0}^\infty a_k x^k=\frac{h(x)}{(1-x)^{n+1}} \label{5.1}
\end{equation}
with $h(x)=h_0+h_1x+\cdots+h_{n+1}x^{n+1}$ a polynomial of degree at most
$n+1$. For example, the generating function of $a_k=k^n$ has this form
with $h(x)$ the Eulerian polynomials of \textbf{(H4$''$)}. Rational
generating functions characterize sequences $\{a_h\}$ which satisfy a
constant coefficient recurrence \cite{Stan}. They arise naturally as the
Hilbert series of graded algebras (\cite{Eis}, Chapter 10.4).

Suppose we are interested in every $r$-th term $\{a_{rk}\},\,0\leq
k<\infty$. It is not hard to see that
\begin{equation*}
\sum_{k=0}^\infty a_{rk}x^k=\frac{h^{<r>}(x)}{(1-x)^{n+1}}
\end{equation*}
for another polynomial $h^{<r>}(x)$ of degree at most $n+1$. Brenti and
Welker \cite{BW} show that the $i$-th coefficient of $h^{<r>}(x)$
satisfies
\begin{equation*}
h_i^{<r>}=\sum_{j=0}^{n+1}C(i,j)h_j
\end{equation*}
with $C$ an $(n+2)\times(n+2)$ matrix with $(i,j)$ entry ($0\leq i,j\leq
n+1$) equal to the number of solutions to $a_1+\cdots+a_{n+1}=ib-j$ where
$0\leq a_l \leq b-1$ are integers. The carries matrix is closely related
to their matrix. Indeed, remove from $C$ the $i=0,n+1$ rows and the
$j=0,n+1$ columns. Let $i'=i-1,\,j'=j-1$. This gives an $n\times n$ matrix
with $i',j'$ entry ($0\leq i',j'\leq n-1$) equal to the number of
solutions to $a_1+\cdots+a_{n+1}=(i'+1)b-(j'+1)$ where $0\leq a_l \leq
b-1$ are integers. Multiplying by $b^{-n}$ and taking transposes gives the
carries matrix for mod $b$ addition of $n$ numbers (see the top of p. 140
of \cite{Holte}). Brenti and Welker develop some properties of the
transformation $C$. We hope some of the facts from the present development
(in particular the central limit theorems satisfied by the coefficients)
will illuminate their algebraic applications.

\subsection{Carries for multiplication}

Consider the process of base $b$ multiplication of a random number (digits
chosen from the uniform distribution on $\{0,1,\cdots,b-1\}$) by a fixed
number $k>0$. We do not require that $k$ is single-digit. Then there is a
natural way to define a carries process, which is best defined by example.
Let $k=26$ and consider multiplying $1423$ by $26$ base 10. The zeroth
carry is defined as $\kappa_0=0$. To compute the first carry, note that
$26 \times 3=78$, so $\kappa_1=7$. Then $\kappa_1+26 \times 2=59$, so
$\kappa_2=5$. Next $\kappa_2+26 \times 4=109$, so $\kappa_3=10$. Finally,
$\kappa_3+26 \times 1=36$, so $\kappa_4=3$.

It is not difficult to see that the above process is a Markov chain on the
state space $\{0,1,\cdots,k-1\}$. For example, if $b=10$ and $k=7$, the
transition matrix is
\[ K(i,j) =
\frac{1}{10} \left[ \begin{array}{c c c c c c c }
2 & 1 & 2 & 1 & 2 & 1 & 1 \\
2 & 1 & 2 & 1 & 1 & 2 & 1 \\
2 & 1 & 1 & 2 & 1 & 2 & 1 \\
1 & 2 & 1 & 2 & 1 & 2 & 1 \\
1 & 2& 1& 2& 1& 1& 2\\
1 & 2 &1 & 1 & 2 & 1 & 2 \\
1  & 1 & 2 & 1 & 2 & 1 & 2
\end{array} \right] \]

The matrix above $K(i,j)$ does not have all eigenvalues real, but the
following properties do hold in general:

\begin{itemize}
\item $K(i,j)$ is doubly stochastic, meaning that
every row and column sums to $1$.
\item $K(i,j)$ is an generalized circulant matrix, meaning that
each column is obtained from the previous column by shifting it downward
by $b \ mod \ k$.
\item Fix $k$ and let $K_a, K_b$ be the base $a,b$ transition
matrices for multiplication by $k$. Then $K_{ab}=K_a K_b$.
\end{itemize}

The first two properties are at the level of undergraduate exercises, and
Chapter 5 of \cite{Davis} is a useful reference for generalized
circulants. The third property holds for the same reason that it does for
Holte's matrix (see the explanation on page 143 of \cite{Holte}).

Since $K$ is doubly stochastic, the carries chain for multiplication has
the uniform distribution on $\{0,1,\cdots,k-1\}$ as its stationary
distribution. Concerning convergence rates, one has the following simple
upper bound for total variation distance.

\begin{prop} \label{fast} Let $K_0^r$ denote the distribution of the carries chain
for multiplication by $k$ base $b$ after $r$ steps, started at the state
$0$. Let $\pi$ denote the uniform distribution on $\{0,1,\cdots,k-1\}$.
Then \[ \frac{1}{2} \sum_{j=0}^{k-1} |K_0^r(j) - \pi(j)| \leq
\frac{k}{2b^r} .\]
\end{prop}

\begin{proof} Observe that \[ K_0^r(j) = \frac{1}{b^r} \left| \{x: jb^r
\leq kx < (j+1)b^r, 0 \leq x < b^r \} \right|.\] The number of integers
$x$ satisfying $\frac{jb^r}{k} \leq x < \frac{(j+1)b^r}{k}$ is between
$\frac{b^r}{k}-1$ and $\frac{b^r}{k}+1$. Hence $|K_0^r(j) - \pi(j)| \leq
\frac{1}{b^r}$, and the result follows by summing over $j$. \end{proof}

Convergence rate lower bounds depend on the number theoretic relation of
$k$ and $b$ in a complicated way. For instance if $k=b$, the process is
exactly random after 1 step.

\subsection{Another amazing matrix}

From one point of view, Holte's amazing matrix exists because there is a
``big'' Markov chain on the symmetric group $S_n$ with eigenvalues
$1,1/b,1/b^2,\cdots$ and a function $T:S_n\to\{0,1,\cdots,n-1\}$ with
image this very same Markov chain. Of course, the interpretation as
``carries'' remains amazing. There are many functions of the basic riffle
shuffling Markov chain which remain Markov chains. Here is a simple one.
Consider repeated shuffling of a deck of $n$ cards using the
Gilbert--Shannon--Reed $b$-shuffles. The position of card labeled ``one''
gives a Markov chain on $\{1,2,\cdots,n\}$. In \cite{ADS} the transition
matrix of this chain is shown to be
\begin{align}
&Q_b(i,j)=\frac{1}{b^n}\times\\
&\sum_{h=1}^b
\sum_{r=\ell}^u\binom{j-1}{r}\binom{n-j}{i-r-1}h^r(b-h)^{j-1-r}
(h-1)^{i-1-r}(b-h+1)^{(n-j)-(i-r-1)}\notag \label{6}
\end{align}
where the inner sum is from $\ell=\max(0,(i+j)-(n+1))$ to
$u=\min(i-1,j-1)$. For example, when $n=2,3$ the matrices are
\begin{gather*}
\frac{1}{2b}\begin{pmatrix}
b+1&b-1\\
b-1&b+1
\end{pmatrix}\\
\frac{1}{6b^2}\begin{pmatrix}
(b+1)(2b+1)&2(b^2-1)&(b-1)(2b-1)\\
2(b^2-1)&2(b^2+2)&2(b^2-1)\\
(b-1)(2b-1)&2(b^2-1)&(b+1)(2b+1)
\end{pmatrix}.
\end{gather*}
The matrix $Q_b$ is shown to satisfy
\begin{itemize}
\item $Q_b$ has eigenvalues $1,1/b,1/b^2,\cdots,1/b^{n-1}$.
\item The eigenvectors of $Q_b$ do not depend on $b$; in particular,
  the stationary distribution is uniform: $\pi(i)=1/n,\,1\leq i\leq
  n$.
\item $Q_aQ_b=Q_{ab}$.
\end{itemize}
We guess that $Q_b$ has other nice properties and appearances.

\section*{Acknowledgments}
We thank Alexi Borodin, Francesco Brenti, Jim Fill
  and Phil Hanlon for real help with this paper. The work of Fulman was
supported by NSF grant DMS-0503901.

\end{document}